\documentclass[12pt]{article}
\usepackage{amssymb}
\usepackage{color}
\usepackage[dvips]{graphicx}
\usepackage[dvips]{epsfig}
\usepackage{amsmath}

\newtheorem{thm}{Theorem}[section]
\newtheorem{defn}{Definition}[section]
\newtheorem{exa}{Example}

\newtheorem{prop}{Proposition}[section]
\newtheorem{cor}[thm]{Corollary}

\newtheorem{rem}{Remark}[section]
\newtheorem{lem}{Lemma}[section]

\def\ben{\begin{enumerate}}
\def\een{\end{enumerate}}
\def\bexa{\begin{exa}}
\def\eexa{\end{exa}}
\def\eproof{\hfill$\rule{2mm}{2mm}$\par}

\def\bqu{\begin{quote}}
\def\equ{\end{quote}}

\def\¯{\hbox{ƒ\it y}}
\def\{\hbox{ƒ\it z}}

\def\ds{\displaystyle}

\def\proof{\noindent{\bf Proof: }}
\def\eproof{\hfill$\rule{2mm}{2mm}$\par}

\def\g{\hbox{$\gamma$}}
\def\d{\hbox{$\delta$}}

\def\i{\hbox{$\mathbf{i}$}}
\def\j{\hbox{$\mathbf{j}$}}
\def\k{\hbox{$\mathbf{k}$}}

\def\h{{\bf h}}
\def\b{{\bf b}}

\def\z{{\bf z}}

\def\be{\begin{equation}}
\def\ee{\end{equation}}
\newcommand{\ba}{\begin{eqnarray}}
\newcommand{\ea}{\end{eqnarray}}
\def\ds{\displaystyle}

\def\R{\hbox{$\mathbb{R}$}}

\def\N{\hbox{$\mathbb{N}$}}
\def\Z{\hbox{$\mathbb{Z}$}}
\def\C{\hbox{$\mathbb{C}$}}
\def\H{\hbox{$\mathbb{H}$}}
\def\O{\hbox{$\mathbb{O}$}}
\def\gcd{\hbox{\rm{gcd}}}

\def\ca{{\cal A}}

\def\m{{\bf m}}

\def\e{{\bf e}}

\def\c{{\gamma}}
\def\d{{\delta}}

\def\b{\hbox{$\beta$}}
\def\a{\hbox{$\alpha$}}

\def\deg{\hbox{deg}}

\begin{document}
\baselineskip=16.5pt

\baselineskip=17.5pt

\title{\sf Polynomials and degrees of maps in real normed algebras}

\author{Takis~Sakkalis\\
Mathematics Laboratory, Agricultural University of Athens, \\
75 Iera Odos, Athens 11855, GREECE}

\date{\today}

\maketitle

\begin{abstract}
\noindent
Let $\cal{A}$ be the algebra of quaternions $\mathbb{H}$ or octonions $\mathbb{O}$. In this manuscript a new proof is given, based on ideas of Cauchy and D' Alembert, of the fact that an ordinary polynomial $f(t) \in {\cal{A}}\, [t]$ has a root in $\cal{A}$. As a consequence,  the Jacobian determinant 
$|J(f)|$ is always non negative in $\cal{A}$.  Moreover, using the idea of the topological degree we show that a regular polynomial $g(t)$ over $\cal{A}$ has also a root in $\cal{A}$. Finally, utilizing multiplication $(*)$ in $\cal{A}$, we prove various results on the topological degree of products of maps. In particular, if $S$ is the unit sphere in $\cal{A}$ and $h_1, h_2: S \to S$ are smooth maps, it is shown that $\hbox{deg} (h_1 * h_2)=\hbox{deg} (h_1) + \hbox{deg} (h_2)$.  
{{\em MCS 2010:} 26B10, 12E15, 11R52}
\end{abstract}

\thispagestyle{empty}

\bigskip
\centerline{e--mail: stp@aua.gr}

\newpage

\setcounter{page}{1}

\thispagestyle{plain}

\section{Introduction} 
\label{sec:intro}

In this section we will briefly state some preliminaries neeeded for this work. We begin with the description of the normed algebras of quaternions $\H$ and octonions $\O$. 

$\Diamond$ {\sf Quaternions:} 
An element $c$ of $\H$ is of the form $c=c_0 + \i c_1 + \j c_2 + \k c_3$, where $c_i \in \R$ and $\i, \j, \k$ are such that $\i^2=\j^2=\k^2=-1$ and $\i \j =-\j \i=\k, \j \k =-\k \j =\i, \k\i =-\i \k =\j$. 
The real part of $c$ is $Re(c)=c_0$ while the imaginary part $Im(c)=\i c_1 + \j c_2 + \k c_3$. The norm of $c$, $|c|=\sqrt{c_0^2+c_1^2+c_2^2+c_3^2}$, its conjugate $\bar{c}=c_0 - \i c_1 - \j c_2 - \k c_3$ while its  inverse is $c^{-1}=\bar{c}	 \cdot |c|^{-2}$, provided that $|c| \ne 0$. $c$ is called an {\em imaginary} unit if $Re(c)=0$ and $|c|=1$, and it has the property $c^2=-1$. In that regard, multiplication in $\H$ is associative but not commutative. 

$\Diamond$ {\sf Octonions:} An octonion $c$ is an element of the form $c=c_0 + \sum_{k=1}^7 \e_k\,c_k $, where $c_0, c_k \in \R$, and $\e_k^2=-1$. To define the algebra structure of octonions, it is enough to specify the multiplication table for the imaginary elements $\e_1, \cdots, \e_7$. For brevity this can be described as follows: write out seven triples of imaginary elements (1) $\e_1, \e_2, \e_3$; (2) $\e_1, \e_4, \e_5$; (3) $\e_1, \e_6, \e_7$; (4) $\e_2, \e_6, \e_4$; (5) $\e_2, \e_5, \e_7$; (6) $ \e_3, \e_4, \e_7$ and (7) $\e_3, \e_5, \e_6$. In each triple, we multiply elements just in the same way as in quaternions. For example, in triple (3) we have: $\e_1\,\e_6=-\e_6\,\e_1=\e_7, \e_6\,\e_7=-\e_7\,\e_6=\e_1, \e_7\,\e_1=-\e_1\,\e_7=\e_6.$\footnote{There are several ways to define multiplication in $\O$ since the vector product of two elements in $\R^7$ is not unique. We opt to choose this way to conform with multiplication in $\H$--as a {\em natural} subset of $\O$--and identify $\e_1=\i, \e_2=\j$ and $\e_3=\k$.} The real part of $c$ is $Re(c)=c_0$ while the imaginary part $Im(c)=\sum_{k=1}^7 \e_k \, c_k $. The norm of $c$, $|c|=\sqrt{c_0^2+c_1^2+\cdots + c_7^2}$, its conjugate $\bar{c}=c_0 - \sum_{k=1}^7 \e_k\, c_k $ while its  inverse is $c^{-1}=\bar{c} \cdot |c|^{-2}$, provided that $|c| \ne 0$. $c$ is called an {\em imaginary} unit if $Re(c)=0$ and $|c|=1$, and it has the property $c^2=-1$. The algebra of octonions is non commutative, non associative, but is is alternative; that is for every $a, b \in \O$, we have $(ab)b=a(bb)$ and $b(ba)=(bb)(a)$. Moreover, for any octonions $a, b$, we have $( \cdots ((a \underbrace{b)b) \cdots b}_{n \;\hbox{\scriptsize{times}}})=ab^n$.

Throughout this note, $\ca$ will stand for  either $\H$ or $\O$ equipped with multiplication $(*)$ as defined above. Also, $\m$ will either be $3$ or $7$ and it shall not be confused with $m$ which might be used as an index.  
An element $c \in \ca$ can also be represented via a real matrix ${\cal C}$ so that $x \, {\cal C}^t=c * x, \; x=(x_0,  x_1, \cdots, x_{\m})$. For example, if  $c=c_0 + \i c_1 + \j c_2 + \k c_3 \in \H,\;$  ${\cal C}$ has the form: 

$${\cal C}=\left[ \begin{array}{rrrrrr} c_0&-c_1&-c_2&-c_3\\c_1&c_0& -c_3&c_2\\c_2&c_3&c_0&-c_1\\c_3&-c_2&c_1&c_0 \end{array} \right] $$
Notice that $|{\cal C}|=|c|^4$; if $c \in \O, \; |{\cal C}|=|c|^8$. The following notation will be needed in the sequel: 

\begin{defn} For any $k \times k, k=4, 8$ real matrix $B$ and $c \in \H, \O$, we define: $cB \equiv{\cal C} B$ and $Bc \equiv B {\cal C}$, respectively. 
\end{defn}

The elements $c_1, c_2 \in \ca$ are called {\em similar}, and denoted by $c_1 \sim c_2$,  if $c_1 \,\eta =\eta \, c_2$ for  a non zero $\eta \in \ca$. Similarity is an equivalence relation and for $c \in \ca$ let us denote by $[c]$ its equivalence class. The following is a useful criterion of similarity, \cite{topu03}: 

\begin{prop} $c_1, c_2 \in \ca$ are similar if and only if $Re(c_1)=Re(c_2)$ and $|Im(c_1)|=|Im(c_2)|$. Furthermore, any $c \in \ca$ is similar to the complex number $\a +\i \b$ with $\a=Re(c)$ and $|\b|=|Im(c)|$.  
\end{prop}

We may identify $\ca$ with $\R^{\m+1}$ via the map $t=(x_0 + \sum_{k=1}^{\m} \e_k \,x_k) 
 \to (x_0, x_1, \cdots, x_{\m})$. Let $f:\ca \to \ca$. In view of this identification, we can also think of $f$ as a map from  $\R^{\m+1} \to \R^{\m+1}$. Indeed, 
 if $f(x_0 + \sum_{k=1}^{\m} \e_k\, x_k)=f_0(x_0, x_1,  \cdots, x_{\m}) + 
 \sum_{k=1}^{\m}  \e_k \,f_k(x_0, x_1,  \cdots, x_{\m})$ we define $f:\R^{\m+1} \to \R^{\m+1}$ by  $f(x_0, x_1, \cdots, x_{\m})=(f_0, f_1, , \cdots, f_{\m})$.  We can also multiply maps in $\ca$. For example, if $\ca=\H$ and $f, g: \H \to \H, f=(f_0, f_1, f_2, f_3), \, g=(g_0, g_1, g_2, g_3)$, we define $f*g=(f_0g_0-f_1g_1-f_2g_2-f_3g_3, f_0g_1+f_1g_0+f_2g_3-f_3g_2, f_0g_2+f_2g_0+f_3g_1-f_1g_3, f_0g_3+f_3g_0+f_1g_2-f_2g_1)$. 

In addition, if $f: \ca \to \ca$, we have 

\begin{defn} \label{jac} The Jacobian $J(f)(c),\; c \in \ca$ is the matrix $\left [ \frac{\partial f_i}{\partial x_j} \right], i,j=0, \cdots, \m$ evaluated at $c$. The determinant of $J(f)$ will be denoted by $|J(f)|$. 
\end{defn}

\section{Polynomials over $\ca$}


Let $n \in \N \cup \{0\}$ and $t, a_i \in \ca$. A ``monomial'' of degree $n$ is defined as  ${\phi}_a^n(t)=a_0ta_1t \cdots ta_n$. A finite sum of monomials of degree $n$ will be denoted by ${\phi}^n(t)$. In the above, special care has to be taken if $\ca =\O$, where parenthesis are needed to be taken into account in the definition of ${\phi}_a^n(t)$. 

\begin{defn} A polynomial $f(t): \ca \to \ca$ of degree $n$ over $\ca$ is a function of the form  

\be \label{poly} 
f(t)=\sum_{k=0}^n {\phi}^k(t) 
\ee 
$f(t)$ shall be called {\em regular} if either $n=0$ or  $\lim_{|t| \to \infty} |{\phi}^n(t)|= \infty$; otherwise $f(t)$ will be called {\em non regular}. Furthermore, a regular polynomial is called {\em ordinary} if it is of the form 
\begin{enumerate} \item $f(t)=a_nt^n + a_{n-1}t^{n-1} + \cdots + a_0$ or 
\item $f(t)=t^na_n + t^{n-1}a_{n-1} + \cdots + a_0, a_i \in \ca$
\end{enumerate} In the former case  $f$ is called {\em left} while in the latter {\em right}. 
\end{defn}
If $c \in \ca$ with $f(c)=0$, $c$ is called a {\em zero} or a {\em root} of $f$. 
For the sake of brevity we shall call in the sequel, unless otherwise stated, an ordinary  polynomial simply {\sf polynomial}. For convenience, we will work with {\em left} polynomials and note that all the results proven hold true for right polynomials as well. To this end we shall give a new proof, based on ideas of Cauchy and D' Alembert \cite{Bou}, of the fact  that every ordinary polynomial $f(t)$ of positive degree has a root in $\ca$, [Theorem \ref{Thm-left}]. The proof that every regular polynomial of positive degree has also a root in $\ca$ will be deferred to the next section.

First we will need the following: 
\begin{rem} Let $n \in \N$ and $a \in \ca$. Then, the equation $t^n-a$ has  a solution in $\ca$.
\end{rem}
\proof Let $\a +\i \b \in \C$ be similar to $a$; that is, $a=\eta (\a +\i \b) {\eta}^{-1}$. Write $\a +\i \b=|\a +\i \b|\,e^{i \theta}$ and let $b={|\a +\i \b|}^{1/n}\,e^{i \phi}$, where $\phi={\theta}/n$. Then, if $\zeta=\eta b {\eta}^{-1}$ we have ${\zeta}^n=\eta b^n {\eta}^{-1}=a$. 
\eproof 
Now we have: 
\begin{lem} \label{lem1} If $f(t)$ is a polynomial of positive degree $n$, then for every $t_0 \in \ca$ with $f(t_0) \ne 0$ and for every $r >0$, there exists  $t \in \ca$ with $|t - t_0| <r$ so that $|f(t)| < |f(t_0)|$.
\end{lem}
\proof We argue by contradiction. Assume then that there exists an $r_0 >0$ so that for each $|t -t_0| \leq r_0$, $\; 0<|f(t_0)| \leq |f(t)|$. We may assume that $t_0 \in \R$; for if not, replace $t$ with $u=t * t_0$ and thus $f(1)=f(1 * t_0)=f(t_0)$.  
We now consider the polynomial $q(t)=\frac{f(t+t_0)}{f(t_0)}$. Then, $q(t)$ has degree $n$ and constant term equal to $1$. Also, observe that $1=q(0) \leq |q(t)|$ for all $|t| \leq r_0$. Furthermore, $q(t)$ is an ordinary polynomial since $t_0 \in \R$ and thus $q(t)=1 + b_kt^k + \cdots + b_nt^n$ with $b_k \ne 0$. Let $\zeta$ be  a solution of the equation $t^k=-\frac{|b_k|}{b_k}$. Note that $|\zeta|=1$. Let $I=\{r\zeta\,|\, 0 <r \leq r_0\}$. For $r\zeta \in I$ we have $|q(r \zeta)| \leq |1 + b_kr^k {\zeta}^k| + |b_{k+1} r^{k+1} {\zeta}^{k+1}| + \cdots + |b_n r^n {\zeta}^n|$. Now we have $|1 + b_kr^k {\zeta}^k|=|1-r^k |b_k|\,|=1-{\rho}^k \,|b_k|$, for some $\rho < r_0$. Thus we get $|q(\rho \, \zeta)| \leq 1-{\rho}^k (|b_k| -|b_{k+1}| {\rho} - \cdots - |b_n| \, {\rho}^{n-k})$. Now for perhaps an even smaller $0< {\rho}_1 < \rho$ we will have $|q(r \, \zeta)\,| <1$ for $0<r < {\rho}_1$, a contradiction to the fact that $|q(t)| \geq 1$ for all $|t| \leq r_0$. \eproof 

\begin{thm} \label{Thm-left} Let $f(t)=a_nt^n + a_{n-1}t^{n-1} + \cdots + a_0$, with $n \geq 1$. Then, $f$ has a root in $\ca$.
\end{thm}
\proof First note that $\lim_{|t| \to \infty} \, |f(t)|=+ \infty$ because $\lim_{|t| \to \infty} \, |a_n t^n|=+ \infty$. Let now $\gamma=\inf \{|f(t)|\, , t \in \ca \}$. Since $\lim_{|t| \to \infty} \, |f(t)|=+ \infty$, there exists an $r > \gamma$ so that $\gamma=\inf \{|f(t)| :  |t| \leq r\}$.
Since the closed ball $B=\{ t \in \ca\,: |t| \leq r\}$ is compact and the function $t \to |f(t)|$ is continuous, there must be a $t_0, \, |t_0| \leq r$ with $\gamma=|f(t_0)|$. Finally observe that $|f(t_0)| \leq |f(t)|$ for $|t| \leq r$. But if $f(t_0) \ne 0$ this contradicts Lemma \ref{lem1}. \eproof

If $g(t)= b_mt^m + b_{m-1}t^{m-1} + \cdots + b_0$  is another polynomial, their product $fg(t)$ is defined in the usual way: 

$$ \label{prod} fg(t)=\sum_{k=0}^{m+n}c_k t^k, \quad \hbox{where $c_k=\ds\sum_{i=0}^k a_i b_{k-i}$} $$
Note that in the above setting the multiplication is performed as if the coefficients  were chosen in a commutative field. However, due to the non commutative nature of $\ca$, we have that $(fg)(t) \ne f(t) * g(t)$, when \be (fg)(t)=F_0 + \sum_{i=1}^{\m} \e_i\,F_i, \;\; f=f_0 + \sum_{i=1}^{\m} \e_i\,f_i, \;\; 
g=g_0 + \sum_{i=1}^{\m} \e_i\,g_i \ee 

According to  Theorem 1 of \cite{gordon64} an element $c\in \ca$ is a zero of $f$ if and only if there exists a polynomial $g(t)$ such that $f(t)=g(t)(t-c)$.  In that way, $f$ can be factored into a product of linear factors $(t-c_i), \; c_i \in \ca$. Indeed, since  $f(t)=g(t)(t-c)$ and $g(t)$ has a root, simple induction shows that 
\be \label{factor} f(t)=a_n(t-c_n)(t-c_{n-1}) \cdots (t-c_1), \; c_j \in \ca \ee

A word of caution: In the above factorization, while $c_1$ is necessarily a root of $f$, $c_j, j=2, \cdots n,$ might not be roots of $f$. For example, the polynomial $f(t)=(t + \k)(t + \j)(t + \i)=t^3+ (\i +\j + \k)t^2+(-\i +\j  - \k)t +1$ has only one root, namely $t=-\i$. Theorem 2.1 of \cite{GS} provides a more detailed version of the above factorization in the case $\ca=\H$.  

Roots of $f$ are distinguished into two types: (i) {\em isolated} and (ii) {\em spherical}. A root $c$ of $f$ is called spherical if and only if its characteristic polynomial $q_c(t)=t^2-2t\, Re(c)+ |c|^2$ divides $f$; for any such polynomial, call $\a_c \pm \i \b_c$ its complex roots. In that case any $\gamma \in \ca$ similar to $c$, is also a root of $f$. For example, if $f(t)=t^2 +1$, any imaginary unit element $c \in \ca$  is a root of $f$. 

\begin{rem} \label{rem1} The polynomial $f(t)$ has a spherical root if and only if it has roots $\a + \i \b, \a - \i \b, \; \a, \b \in \R, \b \ne 0$. 
\end{rem}

If we write $f$ in the form $f(t)=a_0(t) + \sum_{i=1}^{\m} \e_i \, a_i(t)$ we see that $f$ has no spherical roots if and only if $\gcd(a_j)_{j=0}^{\m}=1$. Such an $f$ will be called {\em primitive}. Then, it is easy to see that a primitive $f(t)$ of degree $n$, has at most $n$ distinct roots in $\ca$. 

The conjugate $\bar{f}$ of $f$ is defined as $\bar{f}=a_0(t) - \sum_{i=1}^{\m} \e_i \, a_i(t)$. Note that 
$\bar{f} * f= a_0^2 + a_1^2 + \cdots + a_{\m}^2$, which is a real positive polynomial. Observe that if $\a +\i \b$ is a root of $\bar{f}*f$, then there exists  $c \in \ca$, similar to $\a +\i \b$ so that $f(c)=0$ [~Theorem 4, p. 221 of \cite{gordon64}~].  

\begin{defn} \label{defn1} Let $\phi(t) \in \C[t]$ and $\zeta \in \C$ be a root of $\phi$. We denote by $m(\phi)(\zeta)$ the multiplicity of $\zeta$. Now let $c \in \ca$ be a root of $f$ and let $\mu=m(\bar{f}*f)(\a_c + \i \b_c)$. Then, (1) if  $c$ is isolated, we define its multiplicity $m(f)(c)$, as a root of $f$, to be $\mu$; (2) if $c$ is spherical, its multiplicity is set to be $2\mu$. 
\end{defn}  
Note that: 
\begin{rem} \label{rem2} Let $t_0$ be a root of $f$ and write $f(t)=g(t)(t-t_0)$. Then, $t_0$ is simple (multiple) if and only if $g(t_1) \ne 0, \; (g(t_1) = 0)$ for any $t_1 \sim t_0$, respectively. Moreover, a {\em primitive} polynomial $f$ has simple roots if and only if $\bar{f}* f$ does.  
\end{rem}

\subsection{$|J(f)| \geq 0$ over $\ca$}

In this paragraph we will show that $|J(f)|$ of a polynomial $f$ is non negative over $\O$; the case of $\ca=\H$ is similar, \cite{Sa17}. In particular, we will prove that  if $t_0$ is a root of $f$, $t_0$ is simple if and only if $|J(f)(t_0)| >0$.  Thus, at a multiple root $|J(f)|$ vanishes. 

Let $\deg(f)=n+1, \, n \geq0$. A first indication of  $|J(f)(t)|$ being non negative is when $t \in \R$. Indeed, if $t=r \in \R$, divide $f(t)-f(r)$ by $(t-r)$ to get $f(t) -f(r)=g(t)(t-r)$. Since $r$ commutes with every element of $\ca$, we see that $f(t) -f(r)=g(t)*(t-r)$. Thus, $|J(f)(r)|=|g(r)|^8 \,|I| \geq 0$. 

Now let $t_0={\tau}_0 + \sum_{k=1}^7 e_k{\tau}_k ={\tau}_0 + \tau \in \ca$. Since $J(f)(t_0)=J(f-f(t_0))(t_0)$ we may assume that $f(t_0)=0$. Further, by replacing $f$ with $f(t+{\tau}_0)$ we see that $f({\tau})=0$. Now let $\eta \in \ca, \, |\eta|=1$ so that $\eta* \tau * {\eta}^{-1}=\i s, \, s=|\tau|$. If $u={\eta}^{-1} t {\eta}$ and $F(t)=(f \circ u)(t)$, then $F(\i s)=0$ and $|J(F)|=|J(f(u))|\,|J(u)|$. But $|J(u)|=1$ since $J(u)$ is nothing but an orthogonal matrix. Finally, by replacing $t$ with $t/s$ in $F(t)$ we may assume that $f(\i)=0$. Therefore, it is enough to show that $J(f)(\i)| \geq 0$. 

Divide $f(t)$ by $t-\i$ to get 
$f(t)=g(t)(t-\i)$. Let $g(t)= b_nt^n + \cdots + b_1t +b_0, \; b_k \in \ca$. We write $f(t)=b_0(t-\i) + b_1t(t-\i) + \cdots + b_nt^n(t-\i)$. Let $A$ be the matrix  $$ \label{ma} A= \left [ \begin{array}{rrrrr} N&0&0&0\\ 0&-N&0&0 \\0&0&-N&0\\0&0&0&-N \end{array} \right ], \quad \hbox{where} \quad N= \left [ \begin{array}{rrrrr} 0&-1 \\ 1&0 \end{array} \right ] $$
Notice that $A^2=-I$.  Furthermore, we have 
\begin{lem} \label{lemj}  
$J(t^k(t-\i))(\i)=A^k$ for $k \geq 0$. 
\end{lem}
\proof We use induction on $k$. When $k=0$ we have $J(t- \i)=I=A^0$. Now we have $t^{k+1}(t-\i)=(t^{k+1}-\i t^k)t$. Set $t^{k+1}-\i t^k=a^0 + \sum_{i=1}^7\e_ia^i, \; t=x_0 + \sum_{i=1}^7\e_ix_i$. Then, a calculation shows that $$ J(t^{k+1}(t-\i))(\i)=  
 \left [ \begin{array}{cccccc} M_1&0&0&0\\ 0&M_2&0&0 \\0&0&M_3&0\\0&0&0&M_4 \end{array} \right ]$$
 
where $$ M_1= \left [ \begin{array}{rrrrr} -a_{x_0}^1&-a_{x_1}^1 \\[0.3cm] a_{x_0}^0&a_{x_1}^0 \end{array} \right ],  M_2= \left [ \begin{array}{rrrrr} a_{x_2}^3&a_{x_3}^3 \\[0.3cm] -a_{x_2}^2&-a_{x_3}^2 \end{array} \right ],  M_3= \left [ \begin{array}{rrrrr} a_{x_4}^5&a_{x_5}^5 \\[0.3cm] -a_{x_4}^4&-a_{x_5}^4 \end{array} \right ], M_4= \left [ \begin{array}{rrrrr} a_{x_6}^7&a_{x_7}^7 \\[0.3cm] -a_{x_6}^6&-a_{x_7}^6 \end{array} \right ]
$$
But then, $J(t^{k+1}(t-\i))(\i)=A \,(J(t^k(t-\i))(\i))=A \cdot A^k$. This finishes induction and the proof. \eproof 

In view of the above Lemma we get $$ \label{jd} J(f)(\i)=b_0I+b_1A+b_2A^2 + \cdots + b_nA^n=\sum_{k=0}^n (-1)^kb_{2k}I + \sum_
{l=0}^n (-1)^l b_{2l+1}A. 
$$ 
Set $\sum_{k=0}^n (-1)^kb_{2k}=B_e, \; \sum_
{l=0}^n (-1)^l b_{2l+1}=B_o$. We claim that $|B_eI + B_oA| \geq 0$. Indeed, if either of $B_e, B_o$ is zero there is nothing to prove. Suppose then $B_e*B_o \ne 0$. Then it is enough to show $|I + CA| \geq 0$ for $C=B_o/B_e$. If $C=\g_0 + \sum_{i=1}^7 \e_i{\g}_i=\g_0 + \g 
,\;$ a calculation--via Maple--shows that 
$$|I + CA|=\left [(1-|\g|^2)^2 + 2\g_0^2(1+|C|^2) \right ] \left [(\g_1-1)^2+|C|^2- \g_1^2 \right ]^2 $$

The above proves the claim. Moreover, $|I + CA|$ vanishes precisely 
when $C=\g, \, |\g|=1$ or $C=\i$; that is $C$ is an imaginary unit. In short, $|B_eI + B_oA|=0$ if and only $B_e +B_o\,\d=0$, for a suitable imaginary unit $\d$. 

Let $\d \in \ca$ be an  imaginary unit. Recall that $\d \sim \i$. Then, $g(\d)=B_e + B_o \d$, since $\d^2=-1$. Thus, if $g(\d) \ne 0$, which in turn says that $m(f)(\i)=1$, $|B_eI + B_oA| > 0$. On the other hand, if $g(\g) =0 $, which means that $m(f)(\i) \geq 2$, 
then $|J(f)(\i)|$ vanishes, as required. We summarize the above into the following:
\begin{thm} \label{thm-2}Let $f(t)$ be a polynomial of positive degree. Then, $|J(f)(t)| \geq 0$ for all $t \in \ca$. Moreover, if $t_0$ is a root of $f$, $t_0$ is {\em simple} if and only if $|J(f)(t_0)| > 0$. 
\end{thm}

\section{Degrees of maps in $\ca$} 

Eilenberg and Niven in \cite{EN}, proved the fundamental theorem of algebra for quaternions using a degree argument. The key ingredient was Lemma 2 whose proof was depended on the positiveness of $|J(t^n)|$ at the roots of the equation $t^n=\i$ along with the fact that $\deg(t^n)=n$ for $t \in S^4$. In this section we will prove a slightly more general result, namely that any regular polynomial of positive degree over $\ca$ has a root in $\ca$. In addition, we will show that the topological degree of products of maps is additive.

Let $f: \ca \to \ca$ be a smooth  map with the property that $\lim_{|t| \to \infty} \, |f(t)|=\infty$.  Then, if $\Sigma$ is the  spherical compactification of $\ca$ [$\,\Sigma=S^{\m+1}\,$], $f$ can be smoothly extended to give $\hat{f}: \Sigma \to \Sigma$; i.e. $\hat{f}(t)=f(t), \; t \in \ca$ and  $\hat{f}(\infty)=\infty$. In that case we define $\deg(f):=\deg(\hat{f})$.
Now suppose that $g : \ca \to \ca$ is also smooth  with the same property as $f$. Consider the map $F(t):=f(t) * g (t): \ca \to \ca$. Obviously $\lim_{|t| \to \infty} \, |F(t)|=\infty$. Here is our first result: 
\begin{prop} Let $f,g,F$ be as above. Then, $\deg(F)=\deg(f) + \deg (g)$. 
\end{prop}
\proof According to Sard's Theorem, there exist $\h, \z \in \ca$ so that: (1)  $\h, \z$ are regular values of $f, g$ and (2) the equations $f(t)-\h=0, \, g(t)- \z=0$ have no common roots. Since $\Sigma$ is compact, the sets $A=\{a \in \ca\,|\, f(a)- \h=0\}$ and  $B=\{b \in \ca\,|\, g(b)- \z=0\}$ are finite. Then, $\deg(f)=\sum_{a \in A} sign |J(f)(a)|$ and 
$\deg(g)=\sum_{b \in b} sign |J(g)(b)|$. Let $\Phi(t)=(f(t)-\h) * (g(t)- \z)$. Obviously, $\lim_{|t| \to \infty} \, |\Phi(t)|=\infty$.
In that case for $a \in A, \; b \in B$ we have $$
sign |J(\Phi)(a)|= sign |J(f)(a)| \cdot |g(a)- \z|^{\m+1},\;    
sign |J(\Phi)(b)|= sign |J(g)(b)| \cdot |f(b)- \h|^{\m+1}$$

The above calculation shows that $0$ is a regular value of $\Phi$ and thus $\deg(\Phi)=\deg(f) + \deg(g)$. Finally, since $\Phi=f*g-f* \,\z-\h * \, g + \h * \,\z$, $\Phi$ is homotopic to $F$ via the homotopy $\phi_r(t)= (f*g) (t)+(1-r)((-f* \,\z-\h * \, g + \h * \,\z)(t))$ for $t \in \ca$ and $\phi_r(\infty)=\infty$, $0 \leq r \leq 1$. 
\eproof 

Using the result above and simple induction, we get: 
\begin{cor} \label{cor1} The degree of $t^n: \ca \to \ca$ is equal to $n$.   
\end{cor}

Now, we are ready to prove the Fundamental Theorem of Algebra over $\ca$. 
\begin{thm} Any {\em regular} polynomial of positive degree over $\ca$ has a root in $\ca$. 
\end{thm} 
\proof Let $f(t)=\sum_{k=0}^n {\phi}^k(t)$ be regular with $n \geq 1$. Then, since $\lim_{|t| \to \infty} |{\phi}^n(t)|= \infty$, we see that $\lim_{|t| \to \infty} |f(t)|= \infty$ as well. A slight modification of the proof of Lemma 1 of \cite{EN} shows that $f$ is homotopic to the map $t^n: \ca \to \ca$. Corollary \ref{cor1} shows that $\deg(f) =n$ and therefore $f$ is onto. \eproof
Regularity is a necessary condition for a polynomial to have roots as the following example indicates: 
\begin{exa} The polynomial $f(t)=\i t^2 \j + \j t^2 \i -1$ has no roots over $\H$.
\end{exa}
In the above, observe that $\lim_{|t| \to \infty} |\i t^2 \j + \j t^2 \i|$ does not exist. 

Let now $S$ denote the unit sphere in $\ca$; i.e. $S=S^\m$. Obviously $S$ is equipped with the multiplication ($*$) in $\ca$. Let $k \in \Z$ and consider the map $h(t): S \to S, \; h(t)=t^k$. We then have: 
\begin{lem} The degree of $h$ is equal to $k$. 
\end{lem}
\proof If $k=0$, $h$ is constant and thus has degree $0$. Suppose first that $k \geq 1$. For a fixed $0<r<1$, consider the polynomial $p(t)=t^k -\i r$. Then, $p$ has simple roots [~Remark \ref{rem2}~] ${\rho}_i \in \ca, \; i=1, \cdots, k$ with ${\rho}_i \in Int(S)$ since $|{\rho}_i|=r^{1/k} <1$. Further, from Theorem \ref{thm-2} we get 
$|J(p)({\rho}_i)|>0$ for each $i$. Let $g(t): S \to S$ be defined by $g(t)=\frac{p(t)}{|p(t)|}$. In that case, $\deg(h)=\sum_{i=1}^k \, sign \, |J(p)({\rho}_i)| =k$, \cite{Mi}, Lemma 3, p. 36.
We now claim that $|g(t)-h(t)| < 2$ for $t \in S$. Indeed, note first that $|g(t)-h(t)| \leq 2$. Furthermore, if for some $t_0 \in S$ $\; |g(t_0)-h(t_0)|=2$, we must have $g(t_0)=-h(t_0)$; that is $t_0^k(1+ |t_0^k -\i r|)=\i r$, a contradiction to $0 <r<1$. Thus, $|g(t)-h(t)| < 2$ and this shows that $h$ and $g$ are homotopic, \cite{Mi}, p. $52$. Therefore, $\deg(g)=\deg(h)=k$.

Let now $k=-1$. Then, $h(t)=t^{-1}$ is nothing but a composition of $\m$ reflections. Indeed, if $S=S^3$ and $f_1=(x,-y,z,w), f_2=(x,y,-z,w), f_3=(x,y,z,-w)$, $h(x,y,z,w)=(x,-y,-z,-w)=(f_1 \circ f_2 \circ f_3)(x,y,z,w)$. The case of $S=S^7$ is similar. Therefore, $\deg(h)=(-1)^\m$. 

Finally, let $k \leq -2$. If $f(t)=t^{-1}$ we observe that $(h \circ f)(t)=t^{-k}$ and thus $\deg(h) \, \deg(f)=-k$, or $\deg(h)=k$. This finishes the proof. \eproof 

Let now $f,g: S \to S$ be smooth and $F: S \to S$ be defined by $F(t)=f(t) * g(t)$. Then, 
\begin{thm} $\deg(F)=\deg(f) + \deg(g)$. 
\end{thm}  
\proof Let $\deg(f)=n, \deg(g)=k$ According to Hopf's Theorem, \cite{Mi}, $f, g$ are smoothly homotopic to $f_1, g_1: S \to S$ where $f_1(t)=t^n$ and $g_1(t)=t^k$. Let ${\phi}_r(t), {\psi}_r(t): [0,1| \times S \to S$ be smooth maps so that  ${\phi}_r(0)=f(t), {\phi}_r(1)=f_1$ and ${\psi}_r(0)=g(t), {\psi}_r(1)=g_1$. Define ${\Phi}_r(t)={\phi}_r(t) * {\psi}_r(t): [0,1| \times S \to S$. Notice that ${\Phi}_r(t)$ is continuous and ${\Phi}_r(1)=f(t) * g(t), \; {\Phi}_r(0)=f_1(t) * g_1(t)$. Thus, $f*g$ is homotopic to the map $t \to t^{n+k}$. The latter implies that $\deg(F)=n +k=\deg(f) + \deg(g)$ as required. \eproof


\def\AMM{{\it Amer.\ Math.\ Monthly\ }}
\def\ACM{{\it Adv.\ Comp.\ Math.\ }}
\def\ACMTMS{{\it ACM Trans.\ Math.\ Software\ }}
\def\ACMTOG{{\it ACM Trans.\ Graphics\ }}
\def\BAMS{{\it Bull.\ Amer.\ Math.\ Soc.\ }}
\def\BHMS{{\it Bull.\ Hellenic.\ Math.\ Soc.\ }}
\def\CAD{{\it Comput.\ Aided Design }}
\def\CAEJ{{\it Comput.\ Aided Eng. J.\ }}
\def\CAGD{{\it Comput.\ Aided Geom.\ Design }}
\def\CAVW{{\it Comput.\ Anim.\ Virt.\ Worlds }}
\def\CG{{\it Computers \& Graphics }}
\def\CVGIP{{\it Comput.\ Vision, Graphics, Image\ Proc.\ }}
\def\GM{{\it Graph.\ Models\ }}
\def\IBMJRD{{\it IBM J.\ Res.\ Develop.\ }}
\def\JCAM{{\it J.\ Comput.\ Appl.\ Math.\ }}
\def\JGG{{\it J.\ Geom.\ Graphics}}
\def\JMAA{{\it J.\ Math.\ Anal.\ Appl.\ }}
\def\JSC{{\it J.\ Symb.\ Comput.\ }}
\def\MC{{\it Math.\ Comp.\ }}
\def\MMAS{{\it Math.\ Methods\ Appl.\ Sci.\ }}
\def\NA{{\it Numer.\ Algor.\ }}
\def\PAMS{{\it Proc.\ Amer.\ Math.\ Soc.\ }}
\def\SIAMJNA{{\it SIAM J.\ Numer.\ Anal.\ }}
\def\SIAMR{{\it SIAM Rev.\ }}
\def\JMAS{{\it J.\ Math.\ Sciences.\ }}
\def\GMJ{{\it Georgian \ Math.\ Journal.\ }}
\def\TAMS{{\it Trans.\ Amer.\ Math.\ Soc.\ }}
\def\MJM{{\it Milan\ J.\ Math.\ }}

\begin{flushleft}

\end{flushleft}

\end{document}